\documentclass[reqno]{amsart}
\usepackage{amsthm,amsmath,amssymb}
\usepackage{stmaryrd,mathrsfs}
\usepackage[T1]{fontenc}
\usepackage{esint}
\usepackage{tikz}

\allowdisplaybreaks

\newcommand{\ud}[0]{\,\mathrm{d}}

\newcommand{\abs}[1]{|#1|}
\newcommand{\Babs}[1]{\Big|#1\Big|}
\newcommand{\Norm}[2]{\|#1\|_{#2}}

\newcommand{\BNorm}[2]{\Big\|#1\Big\|_{#2}}
\newcommand{\pair}[2]{\langle #1,#2 \rangle}

\newcommand{\ave}[1]{\langle #1\rangle}

\newcommand{\BMO}[0]{\operatorname{BMO}}

\newcommand{\loc}[0]{\operatorname{loc}}

\newcommand{\R}{\mathbb{R}}
\newcommand{\C}{\mathbb{C}}

\newcommand{\Exp}[0]{\mathbb{E}}
\newcommand{\eps}[0]{\varepsilon}

\swapnumbers \numberwithin{equation}{section}

\theoremstyle{plain}
\newtheorem{theorem}[equation]{Theorem}
\newtheorem{proposition}[equation]{Proposition}

\newtheorem{lemma}[equation]{Lemma}
\newtheorem{conjecture}[equation]{Conjecture}

\theoremstyle{definition}

\theoremstyle{remark}

\makeatletter
\@namedef{subjclassname@2010}{%
  \textup{2010} Mathematics Subject Classification}
\makeatother

\begin{document}

\title{Of commutators and Jacobians}

\author{Tuomas P.~Hyt\"onen}
\address{Department of Mathematics and Statistics, P.O. Box~68 (Pietari Kalmin katu 5), FI-00014 University of Helsinki, Finland}
\email{tuomas.hytonen@helsinki.fi}


\thanks{The author is supported by the Academy of Finland via project Nos.~307333 (Centre of Excellence in Analysis and Dynamics Research) and 314829 (Frontiers of singular integrals).}
\keywords{Commutator, Beurling transform, Jacobian determinant.}
\subjclass[2010]{42B20, 42B25, 42B37, 35F20, 47B47}


\maketitle

\begin{center}
Dedicated to Professor Fulvio Ricci
\end{center}


\begin{abstract}
I discuss the prescribed Jacobian equation $Ju=\det\nabla u=f$ for an unknown vector-function $u$, and the connection of this problem to the boundedness of commutators of multiplication operators with singular integrals in general, and with the Beurling operator in particular. A conjecture of T.~Iwaniec regarding the solvability for general datum $f\in L^p(\R^d)$ remains open, but recent partial results in this direction will be presented. These are based on a complete characterisation of the $L^p$-to-$L^q$ boundedness of commutators, where the regime of exponents $p>q$, unexplored until recently, plays a key role. These results have been proved in general dimension $d\geq 2$ elsewhere, but I will here present a simplified approach to the important special case $d=2$, using a framework suggested by S.~Lindberg.
\end{abstract}

\tableofcontents

\section{The prescribed Jacobian problem}

Given a vector-valued function $u=(u_j)_{j=1}^d\in \dot W^{1,pd}(\R^d)^d$ in the homogeneous Sobolev space
\begin{equation*}
  \dot W^{1,pd}(\R^d)=\{v\in L^1_{\loc}(\R^d): \partial_i v\in L^{pd}(\R^d)\ \forall i\},
\end{equation*}
it is clear that its Jacobian determinant---a linear combination of $d$-fold products of the various $\partial_i u_j$---satisfies $Ju:=\det\nabla u:=\det(\partial_i u_j)_{i,j=1}^d\in L^p(\R^d)$.

Our starting point is the reverse question: Given $f\in L^p(\R^d)$, is there $u\in \dot W^{1,pd}(\R^d)^d$ such that $Ju=f$? This is a nonlinear PDE, known as the ``prescribed Jacobian equation''. It has been mostly studied for {\em smooth} functions $f$ on {\em bounded} domains $\Omega$ \cite{DM,Ye:1994}, in which case there are signifcant geometric applications (e.g. \cite{Avila:2010}). In the global $L^p$ case that we discuss, there is:

\begin{conjecture}[\cite{Iwaniec:Escorial}]\label{conj:Iwaniec}
For $p\in(1,\infty)$, there exists a continuous $E:L^p(\R^d)\to\dot W^{1,pd}(\R^d)^d$ such that $J\circ E=I$. 
\end{conjecture}

As suggested in \cite{Iwaniec:Escorial}, such an $E$ could be interpreted as a ``fundamental solution of the Jacobian equation''.

The case $p=1$ had already been addressed a little earlier. In this case, a simple integration by parts confirms that 
\begin{equation*}
 u\in \dot W^{1,d}(\R^d)^d\quad\Rightarrow\quad\int Ju=0\quad\Rightarrow\quad J(\dot W^{1,d}(\R^d)^d)\subsetneq L^1(\R^d).
\end{equation*}
A somewhat more careful argument yields:

\begin{theorem}[\cite{CLMS}]\label{thn:CLMSbd}
For $u\in \dot W^{1,d}(\R^d)^d$, $d\geq 2$, we have
\begin{equation*}
  \Norm{Ju}{H^1(\R^d)}\lesssim\Norm{u}{\dot W^{1,d}(\R^d)^d}^d
\end{equation*}
where $H^1(\R^d)$ denotes the Hardy space.
\end{theorem}

Again in the reverse direction, \cite{CLMS} asked: Given $f\in H^1(\R^d)$, is there $u\in \dot W^{1,d}(\R^d)^d$ such that $Ju=f$? As a partial positive evidence, they proved:

\begin{theorem}[\cite{CLMS}]\label{thm:CLMSrep}
For every $f\in H^1(\R^d)$, there are $u^i\in \dot W^{1,d}(\R^d)^d$ and $\alpha_i\geq 0$ such that
\begin{equation*}
  f=\sum_{i=1}^\infty\alpha_i J(u^i),\quad \Norm{u^i}{\dot W^{1,d}(\R^d)^d}\leq 1,\quad\sum_{i=1}^\infty\alpha_i\lesssim\Norm{f}{H^1(\R^d)}.
\end{equation*}
\end{theorem}

What about the (perhaps more usual) non-homogeneous Sobolev space
\begin{equation*}
\begin{split}
  W^{1,p}(\R^d) &:=\{v\in { L^p}(\R^d): \nabla v\in L^p(\R^d)^d\}, \\
  \subsetneq \dot W^{1,p}(\R^d) &:=\{v\in { L^1_{\loc}}(\R^d): \nabla v\in L^p(\R^d)^d\}.
\end{split}
\end{equation*}
Given $f\in L^p(\R^d)$ (resp. $H^1(\R^d)$ if $p=1$), could we even hope to find $u\in W^{1,pd}(\R^d)^d$ with $Ju=f$? It was only fairly recently that this was shown to fail, and in fact quite miserably:

\begin{theorem}[\cite{Lindberg:2017}]\label{thm:Lindberg}
The set
\begin{equation*}
  \Big\{\sum_{i=1}^\infty\alpha_i J(u^i):\Norm{u^i}{W^{1,pd}(\R^d)^d}\leq 1,\ \sum_{i=1}^\infty\abs{\alpha_i}<\infty\Big\},
\end{equation*}
which obviously contains the image  $JW^{1,pd}(\R^d)^d$,
has {\em first category} in $L^p(\R^d)$ if $p\in(1,\infty)$ resp. in $H^1(\R^d)$ if $p=1$.
\end{theorem}

Very roughly speaking, the reason for this negative result is the incompatibility of scaling in $W^{1,pd}(\R^d)^d$ on the one hand, and in $L^p(\R^d)$ if $p\in(1,\infty)$ resp. in $H^1(\R^d)$ if $p=1$ on the other hand, but the precise argument is more delicate.

\section{Functional analysis behind the results}

Both the existence (in Theorem \ref{thm:CLMSrep}) and the non-existence (in Theorem \ref{thm:Lindberg}) of the representation $f=\sum\alpha_i J(u^i)$ are based on the following functional analytic lemma from \cite{CLMS} and its elaboration from \cite{Lindberg:2017}:

\begin{lemma}[\cite{CLMS}]
Let $V\subset X$ be a symmetric bounded subset of a Banach space $X$. Then the following are equivalent:
\begin{enumerate}
  \item\label{it:repr} Every $x\in X$ can be written as $x=\sum_{k=1}^\infty\alpha_k v_k$, where $v_k\in V$, $\alpha_k\geq 0$ and $\sum_{k=1}^\infty\alpha_k<\infty$.
  \item\label{it:norming} $V$ is {\em norming} for $X^*$, i.e.,
$\Norm{\lambda}{X^*}\eqsim\sup_{v\in V}\abs{\pair{\lambda}{v}}\quad\forall\lambda\in X^*.$
\end{enumerate}
\end{lemma}
\begin{lemma}[\cite{Lindberg:2017}]
\eqref{it:repr} either holds for all $x\in X$, or in a subset of first category.
\end{lemma}

For the mentioned theorems, these lemmas are applied with the symmetric set $V=J(B)$, where $B=$ unit ball of $\dot W^{1,pd}(\R^d)^d$ or $W^{1,pd}(\R^d)^d$, which is a bounded subset of the Banach space $X= L^p(\R^d)$ or $X=H^1(\R^d)$. Via the equivalent condition \eqref{it:norming}, the well-known dual spaces $X^*=L^{p'}(\R^d)$ or $X^*=\BMO(\R^d)$ enter the considerations.

In order to obtain Theorem \ref{thm:CLMSrep}, \cite{CLMS} proved that

\begin{proposition}[\cite{CLMS}]\label{prop:CLMS}
Let $d\geq 2$. For every $b\in\BMO(\R^d)$, we have
\begin{equation*}
  \Norm{b}{\BMO(\R^d)}\eqsim\sup\Big\{\Babs{\int b J(u)}:\Norm{\nabla u}{d}\leq 1\Big\}.
\end{equation*}
\end{proposition}

The analogous result for $p\in(1,\infty)$ read as follows:

\begin{theorem}[\cite{Hytonen:comm}]\label{thm:HYTrep}
Let $d\geq 2$ and $p\in(1,\infty)$.
For every $f\in L^p(\R^d)$, there are $u^i\in \dot W^{1,dp}(\R^d)^d$ and $\alpha_i\geq 0$ such that
\begin{equation*}
  f=\sum_{i=1}^\infty\alpha_i J(u^i),\quad \Norm{u^i}{\dot W^{1,dp}(\R^d)^d}\leq 1,\quad\sum_{i=1}^\infty\alpha_i\lesssim\Norm{f}{L^p(\R^d)}.
\end{equation*}
\end{theorem}

\begin{proposition}[\cite{Hytonen:comm}]\label{prop:Hyt}
Let $d\geq 2$ and $p\in(1,\infty)$. For every $b\in L^{p'}(\R^d)$, we have
\begin{equation*}
  \Norm{b}{L^{p'}(\R^d)}\eqsim\sup\Big\{\Babs{\int b J(u)}:\Norm{\nabla u}{dp}\leq 1\Big\}.
\end{equation*}
\end{proposition}

\section{Complex reformulation and connection to commutators for $d=2$}

The various results formulated above are valid, as stated, in all dimensions $d\geq 2$, and their proofs in this generality can be found in the quoted references. We now restrict ourselves to dimension $d=2$ in order to discuss an alternative complex-variable approach that is available in this situation, as suggested by \cite{Lindberg:2017}.

For $u=(u_1,u_2)\in\dot W^{1,2p}(\R^2;\R^2)$, let us denote
\begin{equation*}
  h:=u_1+iu_2\in\dot W^{1,2p}(\C;\C),\quad\partial:=\frac12(\partial_1-i\partial_2),\ \bar\partial:=\frac12(\partial_1+i\partial_2).
\end{equation*}
Then, after some algebra, we find that
\begin{equation*}
\begin{split}
  Ju =\det\begin{pmatrix} \partial_1 u_1 & \partial_2 u_1 \\ \partial_1 u_2 & \partial_2 u_2 \end{pmatrix} 
  =\abs{\partial h}^2-\abs{\bar\partial h}^2
  =:\abs{S(v)}^2-\abs{v}^2,
\end{split}
\end{equation*}
where $v:=\bar\partial h\in L^{2p}(\C)$ is in isomorphic correspondence with $h\in \dot W^{1,2p}(\C;\C)$, and  $S$ is the (Ahlfors--)Beurling (or 2D Hilbert) transform
\begin{equation*}
  Sv(z)=-\frac{1}{\pi}\text{p.v.}\int_{\C}\frac{v(y)\ud y_1\ud y_2}{(z-y)^2},
\end{equation*}
which satisfied the fundamental relation $S\circ\bar\partial=\partial$ and maps
$S:L^p(\C)\to L^p(\C)$ bijectively and isometrically for $p=2$ and isomorphically for all $p\in(1,\infty)$.

Let us now see how Proposition \ref{prop:CLMS} and \ref{prop:Hyt} are connected to commutators when $d=2$. By the reformulations just discussed, we have
\begin{equation*}
  \sup\Big\{\Babs{\int b J(u)}:\Norm{u}{\dot W^{1,2p}(\R^2;\R^2)}\leq 1\Big\} 
  \eqsim  \sup\Big\{\Babs{\int b (\abs{Sv}^2-\abs{v}^2)}:\Norm{v}{L^{2p(\C)}}\leq 1\Big\}
\end{equation*}
denoting $v=\bar\partial(u_1+iu_2)$. We claim that the right side can be further written as
\begin{equation}\label{eq:afterPolar}
  \eqsim\sup\Big\{\Babs{\int b (Sv\overline{Sw}-v\overline{w})}:\Norm{v}{L^{2p(\C)}},\Norm{w}{L^{2p(\C)}}\leq 1\Big\}.
\end{equation}
In fact, ``$\leq$'' is obvious, while ``$\gtrsim$'' follows from the elementary polarisation identity
\begin{equation*}
  a\bar b=\frac14\sum_{\eps=\pm 1,\pm i}\eps\abs{a+\eps b}^2,\qquad a,b\in\C,
\end{equation*}
applied pointwise to both $(a,b)=(Sv,Sw)$ and $(a,b)=(v,w)$, which implies that
\begin{equation*}
\begin{split}
  Sv\overline{Sw}-v\overline{w}
  &=\frac14\sum_{\eps=\pm 1,\pm i}\eps\abs{Sv-\eps Sw}^2-\frac14\sum_{\eps=\pm 1,\pm i}\eps\abs{v-\eps w}^2 \\
  &=\frac14\sum_{\eps=\pm 1,\pm i}\eps\Big(\abs{S(v-\eps w)}^2-\abs{v-\eps w}^2\Big),
\end{split}
\end{equation*}
where $\Norm{v-\eps w}{2p}\leq\Norm{v}{2p}+\Norm{w}{2p}\leq 2$ if $\Norm{v}{2p},\Norm{w}{2p}\leq 1$.

Denoting $g:=\overline{Sw}$, we have $\overline{g}=Sw$ and hence $S^*\overline{g}=S^*S w=w$, where we denoted by $S^*$ the conjugate-linear adjoint of $S$ and used the fact that $S^*S$ is the identity. With this substitution, $g\in L^{2p}(\C)$ and $w\in L^{2p}(\C)$ are in isomorphic correspondence, and we have
\begin{equation*}
  \eqref{eq:afterPolar}
  \eqsim\sup\Big\{\Babs{\int b (Sv\cdot g-v\overline{S^*\overline{g}})}:\Norm{v}{L^{2p(\C)}},\Norm{g}{L^{2p(\C)}}\leq 1\Big\}
\end{equation*}
Finally, using the duality $\int\phi\overline{S^*\psi}=\int S\phi\cdot\overline{\psi}$ with $\phi=bv$ and $\psi=\overline{g}$, we have
\begin{equation}\label{eq:enterCommu}
  \int b (Sv\cdot g-v\overline{S^*\overline{g}})
  =\int (b\cdot Sv\cdot g-S(bv)\cdot \overline{\overline{g}})
  =\int g\cdot[b,S]v,
\end{equation}
where we finally introduced the commutator
\begin{equation*}
  [b,S]v=bSv-S(bv).
\end{equation*}

Now the supremum of (the absolute value of) \eqref{eq:enterCommu} over $\Norm{g}{2p}\leq 1$ is the dual norm $\Norm{[b,S]v}{(2p)'}$, and the supremum of this over $\Norm{v}{2p}\leq 1$ is the operator norm
\begin{equation*}
  \Norm{[b,S]}{L^{2p}(\C)\to L^{(2p)'}(\C)}.
\end{equation*}
Summarising the discussion, we have proved:

\begin{lemma}\label{lem:JnormVsCommu}
Let $p\in[1,\infty)$. Then
\begin{equation*}
    \sup\Big\{\Babs{\int b J(u)}:\Norm{u}{\dot W^{1,2p}(\R^2;\R^2)}\leq 1\Big\} 
    \eqsim  \Norm{[b,S]}{L^{2p}(\C)\to L^{(2p)'}(\C)}.
\end{equation*}
\end{lemma}

Thus Propositions \ref{prop:CLMS} and \ref{prop:Hyt}, for $d=2$, are reduced to understanding the norm of the Beurling commutator $[b,S]:L^{2p}(\C)\to L^{(2p)'}(\C)$. When $p=1$, we have $2p=(2p)'=2$, and we are talking about $L^2$-boundedness of commutators, which is a well-studied topic since the pioneering work of \cite{CRW}. When $p\in(1,\infty)$, we have $2p>2>(2p)'$, and we are talking about the boundedness of commutators between different $L^p$ spaces. This, too, has been well studied in the case that the target space exponent is larger (cf. \cite{Janson:1978}), but we are now precisely in the complementary regime. In this case, the result was only achieved very recently.

\section{The commutator theorem}

Complementing various classical results starting with \cite{CRW}, the following result was recently completed in \cite{Hytonen:comm}:

\begin{theorem}\label{thm:commu}
Let $T=S$ with $d=2$, or more generally, let $T$ be any ``uniformly non-degenerate'' Calder\'on--Zygmund operator on $\R^d$, $d\geq 1$.
Let $1<p,q<\infty$ and $b\in L^1_{\loc}(\R^d)$. Then
\begin{equation*}
  [b,T]:L^p(\R^d)\to L^q(\R^d)\quad\text{boundedly}
\end{equation*}
if and only if
\begin{enumerate}
  \item $p=q$ and $b\in\BMO$ \cite{CRW}, or
  \item $p<q\leq p^*$, where $\frac{1}{p^*}:=(\frac{1}{p}-\frac{1}{d})_+$, and $b\in C^{0,\alpha}$ with $\alpha=d(\frac{1}{p}-\frac{1}{q})$, or
  \item $q>p^*$ and $b$ is constant (this and the previous case are due to \cite{Janson:1978}), or
  \item\label{it:p>q} $p>q$ and $b=a+c$, where $c$ is constant and $a\in L^r$ for $\frac1r=\frac1q-\frac1p$ \cite{Hytonen:comm}.
\end{enumerate}
\end{theorem}

Aside from the new regime of exponents $p>q$, another novelty of \cite{Hytonen:comm} (also when $p\leq q$) is the validity of the ``only if'' implication under the fairly general ``uniform non-degeneracy'' assumption on $T$. Recall that \cite{CRW} proved this direction only for the Riesz transfroms, and \cite{Janson:1978,Uchiyama:1978} for ``smooth enough'' kernels, which has been gradually relaxed in subsequent contributions.

The usual Calder\'on--Zygmund size condition requires the upper bound
\begin{equation*}
  \abs{K(x,y)}\leq\frac{c_K}{\abs{x-y}^d}.
\end{equation*}
on the kernel $K$ of $T$.
``Uniform non-degeneracy'' means that we have a matching lower bound essentially over all positions and length-scales, more precisely:
For every $y\in\R^d$ and $r>0$, there is $x$ such that $\abs{x-y}\eqsim r$ and
\begin{equation*}
  \abs{K(x,y)}\geq \frac{c_0}{\abs{x-y}^d}.
\end{equation*}
This is manifestly the case for the Beurling operator, whose kernel $K(x,y)=-\pi^{-1}/(x-y)^2$ satisfies both bounds with an equality.

More generally, Theorem \ref{thm:commu} holds for both
\begin{enumerate}
  \item two-variable kernels $K(x,y)$ (with very little continuity), and
  \item rough homogeneous kernels
\begin{equation*}
  K(x,y)=K(x-y)=\frac{\Omega((x-y)/\abs{x-y})}{\abs{x-y}^d}
\end{equation*}
as soon as $\Omega$ is not identically zero; this was conjectured by \cite{LORR}, who came very close for $p=q$.
\end{enumerate}

We refer the reader to \cite{Hytonen:comm} for the proof of Theorem \ref{thm:commu} in the stated generality; below we give a much simpler argument in the particular case of the Beurling operator $T=S$, which is relevant for the two-dimensional Jacobian problem, as discussed above.

Indeed, for $d=2$, Theorems \ref{thm:CLMSrep} and \ref{thm:HYTrep} are direct corollaries of Theorem \ref{thm:commu} (via the earlier discussion). For $d>2$, they are not direct consequences of Theorem \ref{thm:commu} itself, but they can nevertheless be proved by adapting the ideas of the proof of Theorem \ref{thm:commu}; see again \cite{Hytonen:comm} for details.

\section{The classical implications}\label{sec:CRW}

We begin with a brief discussion of the ``if'' implications in Theorem \ref{thm:commu}:
\begin{enumerate}
  \item The case $p=q$ and $b\in\BMO$ is the only non-trivial ``if'' statement in Theorem \ref{thm:commu}. There are many excellent discussions of this bound (including two entirely different proofs already in \cite{CRW}), so we skip it here.
  \item If $p<q$ and $b\in C^{0,\alpha}$, we only need the size bound $\abs{K(x,y)}\lesssim\abs{x-y}^{-d}$ to see that
\begin{equation*}
\begin{split}
\abs{[b,T]f(x)}
  &=\Babs{\int (b(x)-b(y))K(x,y)f(y)\ud y} \\
  & \leq\int\abs{b(x)-b(y)}\abs{K(x,y)}\abs{f(y)}\ud y \\
  &\lesssim\int\abs{x-y}^\alpha\abs{x-y}^{-d}\abs{f(y)}\ud y.
\end{split}
\end{equation*}
This is a fractional integral with well-known $L^p\to L^q$ bounds!
  \item If $b=c=$ constant, then $[b,T]=0$ is trivially bounded.
  \item If $p>q$ and $b\in L^r$, we use the boundedness of $T:L^p\to L^p$ and $T:L^q\to L^q$ together with H\"older's inequality
\begin{equation*}
  \Norm{bf}{q}\leq\Norm{b}{r}\Norm{b}{p},\qquad\frac{1}{q}=\frac{1}{r}+\frac{1}{p}
\end{equation*}
to see that both $bT$ and $Tb$ individually are $L^p\to L^q$ bounded.
\end{enumerate}

We then turn to the ``only if'' part, starting with the beautiful classical argument of \cite{CRW} for $p=q$. Given a function $b\in L^1_{\loc}(\C)$ and a ball (disc) $B=B(z,r)\subset\C$, we can pick an auxiliary function $\sigma$ with $\abs{\sigma(x)}=1_B(x)$ so that
\begin{equation*}
\begin{split}
  \int_B &\abs{b(x)-\ave{b}_B}\ud x  =\int_B(b(x)-\ave{b}_B)\sigma(x)\ud x,\\
  &=\frac{1}{\abs{B}}\int_B\int_B (b(x)-b(y))\sigma(x)\ud x\ud y \\
  &=\int_B\int_B \frac{b(x)-b(y)}{(x-y)^2}
    \frac{(x-z)^2-2(x-z)(y-z)+(y-z)^2}{\pi r^2}\sigma(x)\ud x\ud y \\
  &=\sum_{i=1}^3\int g_i(x)\Big(\int \frac{b(x)-b(y)}{(x-y)^2}f_i(y)\ud y\Big)\ud x
  =\sum_{i=1}^3\int g_i[b,S]f_i,
\end{split}
\end{equation*}
for suitable functions $f_i,g_i$ with $\abs{f_i(x)}+\abs{g_i(x)}\lesssim 1_B(x)$, whose explicit formulae can be easily deduced from above. Thus
\begin{equation*}
   \int_B\abs{b-\ave{b}_B}
   \leq\sum_{i=1}^3\Norm{[b,S]}{L^p\to L^p}\Norm{f_i}{p}\Norm{g_i}{p'}
   \lesssim\Norm{[b,S]}{L^p\to L^p}\abs{B}^{1/p}\abs{B}^{1/p'}.
\end{equation*}
Dividing by $\abs{B}^{1/p}\abs{B}^{1/p'}=\abs{B}$ and taking the supremum over all $B$ proves that $\Norm{b}{\BMO}\lesssim\Norm{[b,S]}{L^p\to L^p}$.

With a simple modification of the previous display observed by \cite{Janson:1978}, we also find that
\begin{equation*}
  \int_B\abs{b-\ave{b}_B}
   \leq\sum_{i=1}^3\Norm{[b,S]}{L^p\to{L^q}}\Norm{f_i}{p}\Norm{g_i}{{q'}}
   \lesssim\Norm{[b,S]}{L^p\to L^q}\abs{B}^{1/p}\abs{B}^{{1/q'}},
\end{equation*}
where
\begin{equation*}
  \abs{B}^{1/p+1/q'}
  =\abs{B}^{(1/p-1/q)+1}\eqsim \abs{B}\cdot r_B^{d(1/p-1/q)}=\abs{B}\cdot r_B^\alpha.
\end{equation*}
Thus
\begin{equation*}
  \fint_B\abs{b-\ave{b}_B}\lesssim r_B^\alpha,
\end{equation*}
which a well-known characterisation of $b\in C^{0,\alpha}$. For $\alpha>1$, this space has nothing but the constant functions, completing the sketch of the proof of all the classical ``only if'' statements of Theorem \ref{thm:commu}.

\section{The new case $p>q$}

We finally discuss the proof of the ``only if'' implication of Theorem \ref{thm:commu} in the case $p>q$ that was only recently discovered in \cite{Hytonen:comm}. The above estimate
\begin{equation*}
  \int_B\abs{b-\ave{b}_B}\lesssim\abs{B}^{1/p+1/q'}=\abs{B}^{(1/p-1/q)+1}=\abs{B}^{-1/r+1}=\abs{B}^{1/r'}
\end{equation*}
is still true but seems to be useless in this range. How do we even check that a given function is in $L^r$ + constants?

A convenient tool is as follows:

\begin{lemma}[\cite{Hytonen:comm}, Lemma 3.6]\label{lem:DHKY}
If we have the following bound uniformly for cubes $Q\subset\R^d$:
$$\Norm{b-\ave{b}_Q}{L^r(Q)}\leq C,$$
then there is a constant $c\ (=\lim_{Q\to\R^d}\ave{b}_Q)$ such that
\begin{equation*}
  \Norm{b-c}{L^r(\R^d)}\leq C.
\end{equation*}
\end{lemma}

To estimate the local $L^r$ norm, the following result is useful. Depending on one's background, one may like to call it an iterated Calder\'on--Zygmund or atomic decomposition; one can also view it as a toy version of an influential formula of \cite{Lerner:formula}, featuring merely measurable functions in place $L^1(Q_0)$, the median of $b$ in place of the mean $\ave{b}_{Q_0}$, etc. ``Sparse bounds'' of this type have been extensively used in the last few years; the version below is very elementary compared to several recent highlights, but quite sufficient for the present purposes.

\begin{lemma}\label{lem:Lerner}
Given a cube $Q_0\subset\R^d$ and $b\in L^1(Q_0)$, there is a {\em sparse} collection $\mathscr S$ of the family $\mathscr D(Q_0)$ of dyadic subcubes of $Q_0$ such that
\begin{equation*}
  1_{Q_0}(x)\abs{b(x)-\ave{b}_{Q_0}}
  \lesssim\sum_{Q\in\mathscr S}1_Q(x)\fint_Q\abs{b-\ave{b}_Q}.
\end{equation*}
\end{lemma}

A collection of cubes $\mathscr S$ is called {\em sparse} (or almost disjoint) if there are pairwise disjoint {\em major subsets} $E(Q)\subset Q$ for each $Q\in\mathscr S$, meaning that
\begin{equation*}
  E(Q)\cap E(Q')=\varnothing \quad(\forall Q\neq Q'),\qquad\abs{E(Q)}\geq\frac12\abs{Q}.
\end{equation*}

For $L^p$ estimates, sparse is almost as good as disjoint; namely,
\begin{equation}\label{eq:sparseLp}
  \BNorm{\sum_{Q\in\mathscr S}\lambda_Q 1_Q}{p}
  \eqsim\Big(\sum_{Q\in\mathscr S}\lambda_Q^p\abs{Q}\Big)^{1/p},\qquad\forall \lambda_Q\geq 0,
\end{equation}
where equality would hold for a disjoint collection

With these tools at hand, we are ready to prove that $[b,S]:L^p\to L^q$ for $1<q<p<\infty$ only if $b=a+c$, where $a\in L^r$ with $\frac1r=\frac1q-\frac1p$ and $c$ is constant.
For any cube $Q_0\subset\R^d$, we estimate
\begin{equation*}
\begin{split}
  \Norm{b-\ave{b}_{Q_0}}{L^r(Q_0)}
  &\lesssim\BNorm{\sum_{Q\in\mathscr S}1_Q\fint_Q\abs{b-\ave{b}_Q}}{L^r(Q_0)}\qquad\text{(by Lemma \ref{lem:Lerner})} \\
  &\eqsim\Big(\sum_{Q\in\mathscr S}\abs{Q}\Big[\fint_Q\abs{b-\ave{b}_Q}\Big]^r\Big)^{1/r}\qquad\text{(by \eqref{eq:sparseLp})} \\
  &=  \sum_{Q\in\mathscr S}\abs{Q}\lambda_Q\fint_Q\abs{b-\ave{b}_Q}
  =\sum_{Q\in\mathscr S}\lambda_Q\int_Q\abs{b-\ave{b}_Q},
\end{split}
\end{equation*}
with a suitable dualising sequence $\lambda_Q$ such that
\begin{equation}\label{eq:dualSeq}
  \sum_{Q\in\mathscr S}\abs{Q}\lambda_Q^{r'}=1.
\end{equation}
By the same considerations as in Section \ref{sec:CRW} in the case of just one ball $B$, for each of the cubes $Q\in\mathscr S$ above we find functions $f_Q^i$, $g_Q^i$ with
\begin{equation}\label{eq:fiQgiQ}
 \abs{f^i_Q}+\abs{g^i_Q}\lesssim 1_Q
\end{equation}
such that
\begin{equation*}
  \int_Q\abs{b-\ave{b}_Q}= \sum_{i=1}^3 \int g^i_Q[b,S]f^i_Q.
\end{equation*}

Summarising the discussion so far, we have
\begin{equation}\label{eq:keyInterm}
  \Norm{b-\ave{b}_{Q_0}}{L^r(Q_0)}
  \lesssim \sum_{i=1}^3\sum_{Q\in\mathscr S}\lambda_Q \int g^i_Q[b,S]f^i_Q, 
\end{equation}
where the coefficient $\lambda_Q$ and the functions $f_Q^i,g_Q^i$ satisfy \eqref{eq:dualSeq} and \eqref{eq:fiQgiQ}.

We now enter independent random signs $\eps_Q$ on some probability space, and denote by $\Exp$ the expectation. (For the Jacobian theorem in $d>2$: we need to use random $d$th roots of unity at the analogous step, see \cite{Hytonen:comm}.) With the basic orthogonality $\Exp(\eps_Q\eps_{Q'})=\delta_{Q,Q'}$ and H\"older's inequality after observing that
\begin{equation*}
  \frac1r=\frac1q-\frac1p\quad\Rightarrow\quad\frac{1}{r'}=\frac{1}{q'}+\frac{1}{p}\quad\Rightarrow\quad1=\frac{r'}{q'}+\frac{r'}{p},
\end{equation*}
we have
\begin{equation*}
\begin{split}
  RHS\eqref{eq:keyInterm}
  &= \sum_{i=1}^3\Exp\int\Big(\sum_{Q\in\mathscr S}\eps_Q\lambda_Q^{r'/q'} g^i_Q\Big)[b,S]
    \Big(\sum_{Q'\in\mathscr S}\eps_{Q'} \lambda_{Q'}^{r'/p}f^i_{Q'}\Big) \\
  &\lesssim\Norm{[b,S]}{L^p\to L^q}\BNorm{\sum_{Q\in\mathscr S}\lambda_Q^{r'/q'} 1_Q}{q'}
      \BNorm{\sum_{Q\in\mathscr S}\lambda_Q^{r'/p} 1_Q}{p}\qquad\text{(by \eqref{eq:fiQgiQ})} \\
  &\lesssim\Norm{[b,S]}{L^p\to L^q}\Big(\sum_{Q\in\mathscr S}\lambda_Q^{r'} \abs{Q}\Big)^{1/q'}
      \Big(\sum_{Q\in\mathscr S}\lambda_Q^{r'} \abs{Q}\Big)^{1/p}\qquad\text{(by \eqref{eq:sparseLp})} \\
     &=\Norm{[b,S]}{L^p\to L^q}\qquad\text{(by \eqref{eq:dualSeq})}.
\end{split}
\end{equation*}
This shows that
\begin{equation*}
  \Norm{b-\ave{b}_{Q_0}}{L^r(Q_0)}\lesssim \Norm{[b,S]}{L^p\to L^q}
\end{equation*}
for every cube $Q_0$, and hence
\begin{equation*}
  \Norm{b-c}{L^r(\C)}\lesssim \Norm{[b,S]}{L^p\to L^q}
\end{equation*}
for some constant $c$ by Lemma \ref{lem:DHKY}. If we {\em a priori} know that $b\in L^r(\C)$ (as in Proposition \ref{prop:Hyt}), then necessarily $c=0$, and we obtain the desired quantitative bound for $\Norm{b}{L^r(\C)}$.


\begin{thebibliography}{10}

\bibitem{Avila:2010}
A.~Avila.
\newblock On the regularization of conservative maps.
\newblock {\em Acta Math.}, 205(1):5--18, 2010.

\bibitem{CLMS}
R.~Coifman, P.-L. Lions, Y.~Meyer, and S.~Semmes.
\newblock Compensated compactness and {H}ardy spaces.
\newblock {\em J. Math. Pures Appl. (9)}, 72(3):247--286, 1993.

\bibitem{CRW}
R.~R. Coifman, R.~Rochberg, and G.~Weiss.
\newblock Factorization theorems for {H}ardy spaces in several variables.
\newblock {\em Ann. of Math. (2)}, 103(3):611--635, 1976.

\bibitem{DM}
B.~Dacorogna and J.~Moser.
\newblock On a partial differential equation involving the {J}acobian
  determinant.
\newblock {\em Ann. Inst. H. Poincar\'e Anal. Non Lin\'eaire}, 7(1):1--26,
  1990.


\bibitem{Hytonen:comm}
T.~P. Hyt{\"o}nen.
\newblock The {$L^p$-to-$L^q$} boundedness of commutators with applications to
  the {Jacobian} operator.
\newblock {\em Preprint}, arXiv:1804.11167, 2018.

\bibitem{Iwaniec:Escorial}
T.~Iwaniec.
\newblock Nonlinear commutators and {J}acobians.
\newblock {\em J. Fourier Anal. Appl.}, 3:775--796, 1997.

\bibitem{Janson:1978}
S.~Janson.
\newblock Mean oscillation and commutators of singular integral operators.
\newblock {\em Ark. Mat.}, 16(2):263--270, 1978.

\bibitem{Lerner:formula}
A.~K. Lerner.
\newblock A pointwise estimate for the local sharp maximal function with
  applications to singular integrals.
\newblock {\em Bull. Lond. Math. Soc.}, 42(5):843--856, 2010.

\bibitem{LORR}
A.~K. {Lerner}, S.~{Ombrosi}, and I.~P. {Rivera-R{\'{\i}}os}.
\newblock {Commutators of singular integrals revisited}.
\newblock {\em Bull. Lond. Math. Soc.}, 51(1):107--119, 2019.

\bibitem{Lindberg:2017}
S.~Lindberg.
\newblock On the {H}ardy space theory of compensated compactness quantities.
\newblock {\em Arch. Ration. Mech. Anal.}, 224(2):709--742, 2017.

\bibitem{Uchiyama:1978}
A.~Uchiyama.
\newblock On the compactness of operators of {H}ankel type.
\newblock {\em T\^ohoku Math. J. (2)}, 30(1):163--171, 1978.

\bibitem{Ye:1994}
D.~Ye.
\newblock Prescribing the {J}acobian determinant in {S}obolev spaces.
\newblock {\em Ann. Inst. H. Poincar\'e Anal. Non Lin\'eaire}, 11(3):275--296,
  1994.

\end{thebibliography}

\end{document}